\documentclass[]{amsart}

\usepackage{layout}
\usepackage{amsfonts}
\usepackage{amssymb}
\usepackage[mathscr]{eucal}

\evensidemargin\oddsidemargin
\newtheorem{theorem}{Theorem} 

\begin{document}

\title{On compact holomorphically pseudosymmetric \newline K\"ahlerian manifolds}

\author{Zbigniew Olszak}

\keywords{K\"ahlerian manifold, semisymmetry, holomorphic pseudosymmetry}

\begin{abstract}
For compact K\"ahlerian manifolds, the holomorphic pseudosymmetry reduces to the local symmetry if additionally the scalar curvature is constant and the structure function is non-negative. Similarly, the holomorphic Ricci-pseudosymmetry reduces to the Ricci-symmetry under these additional assumptions. We construct examples of non-compact essentially holomorphically pseudosymmetric K\"ahlerian manifolds. These examples show that the compactness assumption cannot be omitted in the above stated theorem. 
\newline
Recently, the first examples of compact, simply connected essentially holomorphically pseudosymmetric K\"ahlerian manifolds are discovered in \cite{J}. In these examples, the structure functions change their signs on the manifold.
\newline
AMS Mathematics Subject Classification (2000) 53C55, 53C25
\end{abstract}

\maketitle

\section{Holomorphic pseudosymmetries}

Let $M$ be a $2n$-dimensional K\"ahlerian manifold with $(J,g)$ as its K\"ahlerian structure. Thus, $J$ is a $(1,1)$-tensor field (an almost complex structure) and $g$ a Riemannian metric on $M$ such that 
$J^2=-I$, $g(J\cdot,J{\cdot\cdot})=g(\cdot,{\cdot\cdot})$ and $\nabla J=0$, $\nabla$ being the Levi-Civita connection of $g$. Let $\mathfrak X(M)$ be the Lie algebra of smooth vector fields on $M$. For $U,V\in\mathfrak X(M)$, let $\mathcal R(U,V) = [\nabla_U,\nabla_V] - \nabla_{[U,V]} = \nabla^2_{UV} - \nabla^2_{VU}$ be the usual curvature operator, and consider additional curvature type operator $\mathcal R^{\mathcal H}(U,V)$ defined by assuming that 
\begin{equation}
\label{tR}
  \mathcal R^{\mathcal H}(U,V)X = g(V,X)U-g(U,X)V+g(JV,X)JU-g(JU,X)JV-2g(JU,V)JX 
\end{equation}
for any $X\in\mathcal X(M)$. The operators $\mathcal R(U,V)$ and $\mathcal R^{\mathcal H}(U,V)$ will be treated as derivations of the tensor algebra on $M$ in the usual sense. For instance, if $T$ is an $(0,k)$-tensor field, then $\mathcal R(U,V)T$, $\mathcal R^{\mathcal H}(U,V)T$ are the $(0,k)$-tensor fields such that 
\begin{eqnarray*}
  (\mathcal R(U,V)T)(X_1,\ldots,X_k)
  &=&-\sum\nolimits_s T(X_1,\ldots,X_{s-1},\mathcal R(U,V)X_s,X_{s+1},\ldots,X_k),\\
  (\mathcal R^{\mathcal H}(U,V)T)(X_1,\ldots,X_k)
  &=&-\sum\nolimits_s T(X_1,\ldots,X_{s-1},\mathcal R^{\mathcal H}(U,V)X_s,X_{s+1},\ldots,X_k).
\end{eqnarray*}
For an $(0,k)$-tensor field $T$, define $(0,k+2)$-tensor fields $\mathcal R\cdot T$, $\mathcal R^{\mathcal H}\cdot T$ by
\begin{eqnarray*}
  (\mathcal R\cdot T)(U,V,X_1,\ldots,X_k)&=&(\mathcal R(U,V)\cdot T)(X_1,\ldots,X_k)\\
  (\mathcal R^{\mathcal H}\cdot T)(U,V,X_1,\ldots,X_k)&=&(\mathcal R^{\mathcal H}(U,V)\cdot T)(X_1,\ldots,X_k).
\end{eqnarray*}

Let us call an $(0,k)$-tensor field $T$ on $M$ to be 
\begin{itemize}
\item
semisymmetric if $\mathcal R\cdot T=0$; 
\item
holomorphically pseudosymmetric if there exists a function $f$ (called the structure function) on $M$ such that $\mathcal R\cdot T = f \mathcal R^{\mathcal H}\cdot T$.
\end{itemize}

A K\"ahlerian manifold will be called 
\begin{itemize}
\item
semisymmetric (resp., Ricci-semisymmetric) if its Riemann (resp., Ricci) curvature tensor is semisymmetric;
\item
holomorphically pseudosymmetric (resp., Ricci-pseudosymmetric) if its Riemann (resp., Ricci) curvature tensor is holomorphically pseudosymmetric. 
\end{itemize}

The class of holomorphically pseudosymmetric K\"ahlerian manifolds contains all semisymmetric K\"ahler\-ian manifolds, especially, those being locally symmetric. For semisymmetric K\"ahlerian manifolds, see among others \cite{BKV,Sin,Sz1,Sz2}. 

The class of holomorphically Ricci-pseudosymmetric K\"ahlerian manifolds contains all Ricci-semi\-symmetric, especially, Ricci-symmetric ($\nabla S=0$), as well as holomorphically pseudosymmetric K\"ahler\-ian manifolds. For Ricci-semisymmetric K\"ahlerian manifolds, see \cite{Mir}.

The holomorphic pseudosymmetry conditions firstly appeared in \cite{ZO}, and after then they were studied in the papers \cite{RD}, \cite{MH}, \cite{J}, \cite{SY}.

It should be said that curvature conditions of this type have also occured under another name in certain papers about projective holomorphic transformations; for some details, see \cite{M,MRH}, etc.

\section{Main results}

Let us start with recalling certain famous examples. Namely, compact 2-di\-men\-sion\-al surfaces, products of compact 2-dimensional surfaces, products of compact 2-dimensional surfaces and complex projective spaces are semisymmetric K\"ahlerian manifolds with non-constant scalar curvature in general. When assuming that they have constant scalar curvatures, they become locally symmetric.

Recently, the problem of the existence of compact essentially  holomorphically pseudosymmetric (that is, different from semisymmetric) K\"ahlerian manifolds was solved in \cite{J}. 

The aim of the presented paper is to prove that under certain additional assumptions, such manifolds do not exist. We also deal with holomorphic Ricci-pseudosymmetry too.

\begin{theorem}
\label{th1}
Let $M$ be a compact K\"ahlerian manifold. Suppose that $M$ is holomorphically Ricci-pseudosymmetric with non-negative structure function $f$, that is, 
\begin{equation}
\label{ricps}
  \mathcal R\cdot S=f \mathcal R^{\mathcal H}\cdot S,\quad f\geqslant0. 
\end{equation}
If the scalar curvature of $M$ is constant, then $M$ is Ricci-symmetric.
\end{theorem}

\begin{theorem}
\label{th2}
Let $M$ be a compact K\"ahlerian manifold. Suppose that $M$ is holomorphically pseudosymmetric with non-negative structure function $f$, that is, 
\begin{equation}
\label{riemps}
  \mathcal R\cdot R=f\mathcal R^{\mathcal H}\cdot R,\quad f\geqslant0. 
\end{equation}
If the scalar curvature of $M$ is constant, then $M$ is locally symmetric.
\end{theorem}

In the last section, we construct examples of holomorphically pseudosymmetric K\"ahlerian manifolds, which are not semisymmetric. For some of them, the scalar curvature is constant and the structure function is positive. This shows that the compactness is an essential assumption in the above theorems. 

\section{Proofs of the theorems}

At first, recall the very well known curvature identities fulfilled by any K\"ahlerian manifold, 
\begin{eqnarray}
\label{beq1}
  & \mathcal R(JU,JV)=\mathcal R(U,V),\quad \mathcal R(JU,V)+\mathcal R(U,JV)=0, &\\
\label{beq2}
  & S(JU,JV)=S(U,V),\quad S(JU,V)+S(U,JV)=0, &\\
\label{beq3}
  & \mathop{\rm Trace}\{X\to \mathcal R(JX,U)V\}=-S(JU,V), &\\
\label{beq4}
  & \mathop{\rm Trace}_g\{(X,Y)\to R(JX,Y,U,V)\}=2S(JU,V), &
\end{eqnarray}
where $R(U,V,X,Y)=g(\mathcal R(U,V)X,Y)$ and $S$ is the Ricci curvature tensor, $S(U,V) = \mathop{\rm Trace}\{X\to \mathcal R(X,U)V\}$. Moreover, the Ricci 2-form $\rho$, $\rho(X,Y)=S(X,JY)$, is closed, and consequently, 
\begin{equation}
\label{ricform}
  (\nabla_XS)(Y,JZ) + (\nabla_YS)(Z,JX) + (\nabla_ZS)(X,JY) = 0.
\end{equation}

\subsection{Proof of Theorem \ref{th1}} 
In our calculations, it will be useful to use the local components tensor convention and the Einstein summation agreement. At first, for the Laplacian of the square of the length of the Ricci tensor $S$, we have 
\begin{equation}
\label{laplas}
 \triangle\big(\|S\|^2\big) = \nabla^i\nabla_i(S_{jk}S^{jk})
   =2(\nabla^i\nabla_iS_{jk})S^{jk} + 2(\nabla_iS_{jk})(\nabla^iS^{jk}).
\end{equation}
In the sequel, we need the following formula
\begin{equation}
\label{xxx}
  -(\nabla_XS)(Y,Z) + (\nabla_YS)(Z,X) + (\nabla_{JZ}S)(X,JY) = 0,
\end{equation}
which can be obtained from (\ref{ricform}) by replacing $Z$ with $JZ$ and next using (\ref{beq2}). In local coordinates, (\ref{xxx}) reads
$$
  \null-\nabla_iS_{jk}+\nabla_jS_{ki}+\nabla_bS_{ia}J^b_kJ^a_j=0.
$$
The covariant differentiation of the above equality gives
$$
  \null-\nabla_h\nabla_iS_{jk}+\nabla_h\nabla_jS_{ki}+\nabla_h\nabla_bS_{ia}J^b_kJ^a_j=0.
$$
Transvecting the last relation with $S^{jk}=g^{ja}g^{kb}S_{ab}$ and using formula $S^{kj}J^b_kJ^a_j=S^{ab}$ (which is a consequence of (\ref{beq2})), we find $(\nabla_h\nabla_iS_{jk})S^{jk}=2(\nabla_h\nabla_jS_{ki})S^{jk}$ and next 
\begin{equation}
\label{yyy}
  (\nabla^i\nabla_iS_{jk})S^{jk}=2g^{hi}(\nabla_h\nabla_jS_{ki})S^{jk}. 
\end{equation}
We are going to transform (\ref{yyy}) by applying the holomorphic Ricci-pseudosymmetry (\ref{ricps}). Using (\ref{tR}) and (\ref{beq2}), we find for $\mathcal R^{\mathcal H}\cdot S$, 
\begin{eqnarray}
\label{tRS}
  (\mathcal R^{\mathcal H}\cdot S)(U,V,X,Y) 
   &=&\null-S(\mathcal R^{\mathcal H}(U,V)X,Y)-S(X,\mathcal R^{\mathcal H}(U,V)Y) \nonumber\\
   &=& \null-g(V,X)S(U,Y)+g(U,X)S(V,Y)-g(V,Y)S(X,U) \nonumber\\
   & &\null+g(U,Y)S(X,V)-g(JV,X)S(JU,Y)+g(JU,X)S(JV,Y) \nonumber\\
   & &\null-g(JV,Y)S(X,JU)+g(JU,Y)S(X,JV).
\end{eqnarray}
Moreover, we have for $\mathcal R\cdot S$, 
\begin{equation}
\label{RS}
  (\mathcal R\cdot S)(U,V,X,Y) 
    = (\mathcal R(U,V)S)(X,Y) = \big((\nabla^2_{UV}-\nabla^2_{VU})S\big)(X,Y) 
    = (\nabla^2_{UV}S)(X,Y)-(\nabla^2_{VU}S)(X,Y). 
\end{equation}
Now, using (\ref{tRS}) and (\ref{RS}) and (\ref{ricps}), we obtain 
\begin{eqnarray}
\label{bum1}
  \nabla_h\nabla_jS_{ki} - \nabla_j\nabla_hS_{ki}
    &=& f\big(\null-g_{jk}S_{hi}+g_{hk}S_{ji}-g_{ji}S_{kh}+g_{hi}S_{kj} \nonumber\\
    &&\quad \null-J_{jk}J^a_h S_{ai}+J_{hk}J^a_jS_{ai}-J_{ji}S_{ka}J^a_h+J_{hi}S_{ka}J^a_j\big), 
\end{eqnarray}
where $J_{ij}=J^a_ig_{aj}(=-J_{ji})$. Note that by (\ref{beq2}), we have 
\begin{equation}
\label{beq2x}
  S_{ab}J^a_iJ^b_j=S_{ij},\quad S_{ia}J^a_j+S_{ja}J^a_i=0.
\end{equation}
From (\ref{bum1}), by tranvection with $g^{hi}$ and using (\ref{beq2x}), it follows that 
\begin{equation}
\label{zzz}
  g^{hi}(\nabla_h\nabla_jS_{ki}) - g^{hi}(\nabla_j\nabla_hS_{ki}) = f(2nS_{jk}-rg_{jk}), 
\end{equation}
where $r$ is the scalar curvature. Since $r$ is constant, it holds $g^{hi}\nabla_hS_{ki}=(1/2)\nabla_kr=0$, and therefore $g^{hi}(\nabla_j\nabla_hS_{ki})=0$. Thus, (\ref{zzz}) leads to 
$$
  g^{hi}(\nabla_h\nabla_jS_{ki}) = f(2nS_{jk}-rg_{jk}),  
$$
which applied to the right hand side of (\ref{yyy}) yields 
$$
  (\nabla^i\nabla_iS_{jk})S^{jk}=2f(2nS_{jk}-rg_{jk})S^{jk}=4nf\big(\|S\|^2-r^2/(2n)\big).
$$
The last equality turns (\ref{laplas}) into 
\begin{equation}
\label{laplas2}
  \triangle\big(\|S\|^2\big) = 8nf\big(\|S\|^2-r^2/(2n)\big) + 2\|\nabla S\|^2. 
\end{equation}

Recall the famous Hopf Lemma, which states that for a function $\varphi$ on a compact Riemannian manifold, if $\triangle\varphi\geqslant0$, then $\triangle\varphi=0$ and the function is constant (cf.\ e.g.\ \cite{KN} or \cite{YB}). 

Returning to our proof, note that for any Riemannian manifold, it always holds $\|S\|^2-r^2/(2n)\geqslant0$. Therefore and by the assumption $f\geqslant0$, the right hand side of (\ref{laplas2}) is non-negative. Consequently, $\triangle\big(\|S\|^2\big)\geqslant0$, and by the Hopf Lemma, $\triangle\big(\|S\|^2\big)=0$. This applied into (\ref{laplas2}) leads to 
$$
  8nf\big(\|S\|^2-r^2/(2n)\big) + 2\|\nabla S\|^2=0.
$$
Hence, it follows that $\|\nabla S\|=0$, and finally $\nabla S=0$, which is just the Ricci-symmetry. This completes the proof.

\subsection{\bf Proof of Theorem \ref{th2}} 
Let $M$ be a holomorphically pseudosymmetric K\"ahlerian manifold with constant scalar curvature and $f\geqslant0$. Since the formula (\ref{riemps}) always implies the condition (\ref{ricps}) with the same structure function, $M$ is holomorphically Ricci-pseudosymmetric. Consequently, by Theorem \ref{th1}, $M$ is Ricci-symmetric, that is , $\nabla S=0$.

To prove that $M$ is in fact locally symmetric, we will use the Lichnerowicz formula, which is valid for any Riemannian manifold (\cite[Lemma 4.7]{Sz1}; see also \cite{AL})
\begin{equation}
\label{lichner}
  \nabla^p\nabla_p\big(R_{ijkl}R^{ijkl}\big) = 2\nabla_pR_{ijkl}\nabla^pR^{ijkl}
  +4R^{ijkl}(\nabla_j\nabla_kS_{il}-\nabla_j\nabla_lS_{ik})-4R^{ijkl}g^{pq}F_{pijqkl},
\end{equation}
where
$$
  F_{pqijkl}=\nabla_p\nabla_qR_{ijkl}-\nabla_q\nabla_pR_{ijkl}. 
$$

For $\mathcal R\cdot R$, we have 
\begin{equation}
\label{RR}
  (\mathcal R\cdot R)(U,V,W,X,Y,Z)=(\mathcal R(U,V)R)(W,X,Y,Z) 
    =((\nabla^2_{UV}R)(W,X,Y,Z)-(\nabla^2_{VU}R)(W,X,Y,Z).
\end{equation}
On the other hand, using (\ref{tR}) and (\ref{beq1}), we find for $\mathcal R^{\mathcal H}\cdot R$, 
\begin{eqnarray}
\label{tRR}
  (\mathcal R^{\mathcal H}\cdot R)(U,V,W,X,Y,Z) 
   &=& (\mathcal R^{\mathcal H}(U,V)R)(W,X,Y,Z) \nonumber\\
  &=&\null-R(\mathcal R^{\mathcal H}(U,V)W,X,Y,Z)-R(W,\mathcal R^{\mathcal H}(U,V)X,Y,Z) \nonumber\\
  &&\quad \null-R(W,X,\mathcal R^{\mathcal H}(U,V)Y,Z)-R(W,X,Y,\mathcal R^{\mathcal H}(U,V)Z). \nonumber\\
  &=&\null-g(V,W)R(U,X,Y,Z)+g(U,W)R(V,X,Y,Z) \nonumber\\
  &&\quad \null-g(JV,W)R(JU,X,Y,Z)+g(JU,W)R(JV,X,Y,Z) \nonumber\\
  &&\quad \null-g(V,X)R(W,U,Y,Z)+g(U,X)R(W,V,Y,Z) \nonumber\\
  &&\quad \null-g(JV,X)R(W,JU,Y,Z)+g(JU,X)R(W,JV,Y,Z) \nonumber\\
  &&\quad \null-g(V,Y)R(W,X,U,Z)+g(U,Y)R(W,X,V,Z) \nonumber\\
  &&\quad \null-g(JV,Y)R(W,X,JU,Z)+g(JU,Y)R(W,X,JV,Z) \nonumber\\
  &&\quad \null-g(V,Z)R(W,X,U,Z)+g(U,Z)R(W,X,V,Z) \nonumber\\
  &&\quad \null-g(JV,Z)R(W,X,Y,JU)+g(JU,Z)R(W,X,Y,JV) 
\end{eqnarray}
Applying (\ref{RR}), (\ref{tRR}) and (\ref{riemps}), we obtain 
\begin{eqnarray*}
  F_{pqijkl}&=&\nabla_p\nabla_qR_{ijkl}-\nabla_q\nabla_pR_{ijkl} \\
            &=& f(\null-g_{qi}R_{pjkl}+g_{pi}R_{qjkl}-J_{qi}J^a_pR_{ajkl}+J_{pi}J^a_qR_{ajkl} \\
            & &\null-g_{qj}R_{ipkl}+g_{pj}R_{iqkl}-J_{qj}J^a_pR_{iakl}+J_{pj}J^a_qR_{iakl} \\
            & &\null-g_{qk}R_{ijpl}+g_{pk}R_{ijql}-J_{qk}J^a_pR_{ijal}+J_{pk}J^a_qR_{ijal} \\
            & &\null-g_{ql}R_{ijkp}+g_{pl}R_{ijkq}-J_{ql}J^a_pR_{ijka}+J_{pl}J^a_qR_{ijka}).
\end{eqnarray*}
From the above, by transvection with $g^{pj}$, we get 
\begin{eqnarray}
\label{uuu}
  g^{pj}F_{pqijkl} 
  &=& f\big((2n-1)R_{iqkl}+R_{ikql}+R_{ilkq}+J^a_iJ^b_qR_{abkl} 
      -J^a_kJ^b_qR_{iabl}-J^a_lJ^b_qR_{iakb} \nonumber\\
  & &\null-g_{qk}S_{il}+g_{ql}S_{ik}
         -J_{qk}J^{ab}R_{ailb}+J_{ql}J^{ab}R_{aikb}+J_{qi}J^{ab}R_{abkl}\big), 
\end{eqnarray}
where $J^{ij}=g^{ia}J^j_a(=-J^{ji})$. We need the following formulas 
\begin{equation}
\label{bumbum}
  J^a_iJ^b_jR_{abkl}=R_{ijkl},\quad J^a_iR_{ajkl}=J^a_jR_{aikl},\quad
  J^{ab}R_{ajkb}=J^a_jS_{ak}, \quad J^{ab}R_{abkl}=-2J^a_kS_{al},
\end{equation}
which are consequences of (\ref{beq1}), (\ref{beq3}) and (\ref{beq4}). Moreover, using the first Binchi identity and (\ref{bumbum}), we can find 
\begin{equation}
\label{bambam}
  R_{ikql}+R_{ilkq} = R_{iqkl},\quad -J^a_kJ^b_qR_{iabl}-J^a_lJ^b_qR_{iakb} = R_{iqkl}.
\end{equation}
By applying (\ref{bumbum}), (\ref{bambam}) and (\ref{beq2x}), we transform (\ref{uuu}) into the following form
\begin{equation}
\label{gF}
  g^{pj}F_{pqijkl} 
  = f\big(2(n+1)R_{kliq}-S_{li}g_{kq}+S_{ki}g_{lq} 
    -J^a_lS_{ai}J_{kq}+J^a_kS_{ai}J_{lq}+2J^a_kS_{al}J_{iq}\big).
\end{equation}
Recall that the holomorphic projective curvature $(1,3)$-tensor $P$ is defined by 
$$
  \mathcal P(U,V)W = \mathcal R(U,V)W-\frac1{2(n+1)}(S(V,W)U-S(U,W)V
         +S(JV,W)JU-S(JU,W)JV-2S(JU,V)JW).
$$
The local coordinates of the $(0,4)$-tensor $P$, $P(W,X,Y,Z)=g(\mathcal P(W,X)Y,Z)$, are the following 
$$
  P_{hijk}=R_{hijk}-\frac1{2(n+1)}(S_{ij}g_{hk}-S_{hj}g_{ik}
           +J^a_iS_{aj}J_{hk}-J^b_hS_{bj}J_{ik}-2J^a_hS_{ai}J_{jk}).
$$
In this context, (\ref{gF}) can be rewritten as
$$
  g^{pj}F_{pqijkl} = 2(n+1)fP_{kliq}
$$
By virtue of the last formula, we obtain
\begin{equation}
\label{RgF}
  R^{ijkl}g^{pq}F_{pijqkl}=2(n+1)fP_{klji}R^{ijkl}=-2(n+1)fP_{lkji}R^{lkji}.
\end{equation}
By straightforward calculations in which (\ref{bumbum}) should be used, we get 
$$
  \|P\|^2=P_{ijkl}P^{ijkl}=P_{ijkl}R^{ijkl}=\|R\|^2-\frac{4}{n+1}\|S\|^2\geqslant0.
$$
Therefore, (\ref{RgF}) can be rewritten as
\begin{equation}
\label{RgF2}
  R^{ijkl}g^{pq}F_{pijqkl}=-2(n+1)f\|P\|^2.
\end{equation}

The already proved condition $\nabla S=0$ and the formula (\ref{RgF2}) enables us to rewrite the Lichnerowicz formula (\ref{lichner}) in the following form 
\begin{equation}
\label{laplas3}
  \triangle\big(\|R\|^2\big)=2\|\nabla R\|^2+8(n+1)f\|P\|^2.  
\end{equation}
As in the previous proof, we use the Hopf Lemma. By the assumption $f\geqslant0$, the right hand side of (\ref{laplas3}) is non-negative. Consequently, $\triangle\big(\|R\|^2\big)\geqslant0$, and by the Hopf Lemma, $\triangle\big(\|R\|^2\big)=0$. This applied into (\ref{laplas3}) leads to 
$$
  2\|\nabla R\|^2+8(n+1)f\|P\|^2=0.
$$
Hence, it follows that $\|\nabla R\|=0$, and finally $\nabla R=0$, which completes the proof.

\section{A class of examples}

Below, we construct a class of examples of non-compact essentially holomorphically pseudosymmetric K\"ahlerian manifolds. For some of them, the scalar curvature is constant and the structure function is positive.

Let $(x^{\alpha},y^{\alpha},z,t)$ denote the Cartesian coordinates in $\mathbb R^{2m+2}$, $m\geqslant1$. Latin indices take on values from 1 to $2m{+}2$, Greek indices will run from 1 to $m$, and $\alpha' =\alpha+m$ for any $\alpha\in\{1,\ldots,m\}$. Assume that $M=N\times(A,B)\subset\mathbb R^{2m+2}$, where $N$ is an open connected subset of $\mathbb R^{2m+1}$, $(A,B)$ is an open interval and $B>A>0$. Suppose that $h\colon(A,B)\to\mathbb R$ is a smooth function which non-zero at any $t\in(A,B)$. Let $(e_i)$ be the frame of vector fields on $M$ defined by 
$$
  e_{\alpha} = \dfrac1t\, \dfrac{\partial}{\partial x^{\alpha}},\ \  
  e_{\alpha'} = \dfrac1t \Big(\dfrac{\partial}{\partial y^{\alpha}}
          +2x^{\alpha}\dfrac{\partial}{\partial z}\Big),\ \ 
  e_{2m+1} = \dfrac1{t^2h}\, \dfrac{\partial}{\partial z},\ \ 
  e_{2m+2} = th\, \dfrac{\partial}{\partial t},
$$
and let $(\theta^i)$ be the dual frame of differential 1-forms, 
$$
  \theta^{\alpha} = t\,dx^{\alpha},\ \ 
  \theta^{\alpha'} = t\,dy^{\alpha},\ \ 
  \theta^{2m+1} = t^2h 
       \Big(\null-2\sum\nolimits_{\lambda} x^{\lambda}dy^{\lambda}+dz\Big),\ \ 
  \theta^{2m+2} = \dfrac{1}{th}\, dt. 
$$
For the non-zero Lie brackets of $e_i$, we have 
$$
  [e_{\alpha},e_{\beta'}] = 2h\delta_{\alpha\beta}\,e_{2m+1}, \quad
  [e_{\alpha},e_{2m+2}] = he_{\alpha}, \quad
  [e_{\alpha'},e_{2m+2}] = he_{\alpha'}, \quad
  [e_{2m+1},e_{2m+2}] = (2h+th')\,e_{2m+1}.
$$
Define an almost complex structure $J$ on $M$ by assuming 
$$
  Je_{\alpha} = e_{\alpha'},\quad Je_{\alpha'}=-e_{\alpha},\quad
  Je_{2m+1} = e_{2m+2},\quad Je_{2m+2}=-e_{2m+1}.
$$
For the Nijenhuis tensor $N_J$, it can be checked that 
$$
  N_J(e_i,e_j)=[Je_i,Je_j]-J[e_i,Je_j]-J[Je_i,e_j]+J^2[e_i,e_j]=0,
$$
for any $i,j$. By the Newlander-Nirenberg theorem, $J$ is a complex structure on $M$. Let $g$ be the Riemannian metric on $M$ for which $(e_i)$ is an orthonormal frame, so that $g = \sum_i \theta^i\otimes\theta^i$. It is obvious that the pair $(J,g)$ is a Hermitian structure on $M$. For the fundamental form $\varOmega$, ${\varOmega}(X,Y)=g(JX,Y)$, we have
\begin{eqnarray*}
  {\varOmega} &=& 2\sum\nolimits_{\lambda}\theta^{\lambda}\wedge\theta^{\widetilde\lambda}
              +2\,\theta^{2m+1}\wedge\theta^{2m+2},\\
         &=& 2t^2\sum\nolimits_{\lambda}dx^{\lambda}\wedge dy^{\lambda}
              +2t\Big(\null-2\sum\nolimits_{\lambda} x_{\lambda}dy_{\lambda}\wedge dt+dz\wedge dt\Big).
\end{eqnarray*}
Hence $d\varOmega=0$, i.e., $\varOmega$ is closed. Thus, the pair $(J,g)$ becomes a K\"ahlerian structure on $M$. We are going to show that it is holomorphically pseudosymmetric. 

For the Levi-Civita connection corresponding to $g$, we have 
\begin{eqnarray*}
  && \nabla_{e_{\alpha}}e_{\beta} = \nabla_{e_{\alpha'}}e_{\beta'} 
      = -h\delta_{\alpha\beta}\,e_{2m+2}, \\
  && \nabla_{e_{\alpha}}e_{\beta'} = -\nabla_{e_{\alpha'}}e_{\beta} 
      = h\delta_{\alpha\beta}\,e_{2m+1}, \\
  && \nabla_{e_{\alpha}}e_{2m+1} = \nabla_{e_{2m+1}}e_{\alpha} 
      = -\nabla_{e_{\alpha'}}e_{2m+2} = -he_{\alpha'}, \\
  && \nabla_{e_{\alpha}}e_{2m+2} = \nabla_{e_{\alpha'}}e_{2m+1} 
      = \nabla_{e_{2m+1}}e_{\alpha'} = he_{\alpha}, \\
  && \nabla_{e_{2m+1}}e_{2m+1} = -(2h+th')e_{2m+2}, \\
  && \nabla_{e_{2m+1}}e_{2m+2} = (2h+th')e_{2m+1}.
\end{eqnarray*}

Let $R_{hijk}(={R_{hij}}^k)$ be the components of the curvature tensor $R$ with respect to the adapted frame, $R(e_h,e_i)e_j=\sum_kR_{hijk}e_k$. The non-zero components of $R$ are related to the following
\begin{eqnarray*}
  && R_{\alpha\beta\gamma\delta}=R_{\alpha\beta\gamma'\delta'}=
       R_{\alpha'\beta'\gamma'\delta'}= h^2(\delta_{\alpha\gamma}\delta_{\beta\delta}-
       \delta_{\alpha\delta}\delta_{\beta\gamma}),\\
  && R_{\alpha\beta'\gamma\delta'}= h^2(\delta_{\alpha\gamma}\delta_{\beta\delta}
       +\delta_{\beta\gamma}\delta_{\alpha\delta}+2\delta_{\alpha\beta}\delta_{\gamma\delta}),\\
  && R_{\alpha\beta'(2m+1)(2m+2)}= 2h(h+th')\delta_{\alpha\beta},\\
  && R_{\alpha(2m+1)\beta(2m+1)}=R_{\alpha(2m+1)\beta'(2m+2)}\\
  &&  \qquad\qquad=R_{\alpha(2m+2)\beta(2m+2)}=-R_{\alpha(2m+2)\beta'(2m+1)}\\
  &&  \qquad\qquad=R_{\alpha'(2m+1)\beta'(2m+1)}=R_{\alpha'(2m+2)\beta'(2m+2)} 
       = h(h+th')\delta_{\alpha\beta},\\
  && R_{(2m+1)(2m+2)(2m+1)(2m+2)}= 4h^2+7thh'+t^2(h'^{\,2}+hh'').
\end{eqnarray*}

On the other hand, for the components of the tensor $\mathcal R^{\mathcal H}$, we have
$$
  \mathcal R^{\mathcal H}_{hijk}\Big(=\mathcal R^{\mathcal H}_{hij}{}^k\Big)
    =g_{hk}g_{ij} - g_{hj}g_{ik} + J_{hk} J_{ij} - J_{hj} J_{ik} - 2 J_{hi} J_{jk},
$$
where $g_{ij}$ and $J_{ij}$ are the components of $g$ and $J$ with respect to $(e_i)$. Thus, $g_{ij}=g(e_i,e_j)=\delta_{ij}$ and $Je_i=\sum_s J_{is}e_s$ with $J_{\alpha\beta'}=-J_{\alpha'\beta}=\delta_{\alpha\beta}$, $J_{(2m+1)(2m+2)}=-J_{(2m+2)(2m+1)}=1$, otherwise $J_{ij}=0$.

The structure $(J,g)$ satisfies the holomorphic pseudosymmetry condition 
$$
  R\cdot R=f\mathcal R^{\mathcal H}\cdot R\quad\mbox{with}\quad f=-h(h+th').
$$
For, it is sufficient to verify that the relation holds $Q\cdot R=0$, where the curvature like tensor $Q=R-f\mathcal R^{\mathcal H}$ is treated as the derivation of the tensor algebra. At first, we find the components of $Q$ with respect to $(e_i)$, which are as follows 
\begin{eqnarray*}
  && Q_{\alpha\beta\gamma\delta}=Q_{\alpha\beta\gamma'\delta'}=
       Q_{\alpha'\beta'\gamma'\delta'}= -thh'
       (\delta_{\alpha\gamma}\delta_{\beta\delta}-\delta_{\alpha\delta}\delta_{\beta\gamma}),\\
  && Q_{\alpha\beta'\gamma\delta'}= -thh'(\delta_{\alpha\gamma}\delta_{\beta\delta}
        +\delta_{\beta\gamma}\delta_{\alpha\delta}+2\delta_{\alpha\beta}\delta_{\gamma\delta}),\\
  && Q_{(2m+1)(2m+2)(2m+1)(2m+2)}= 3thh'+t^2(h'^{\,2}+hh'').
\end{eqnarray*}
Next, we check that the all components 
$$
  (Q\cdot R)_{pghijk} = -\sum\nolimits_s \left(Q_{pqhs} R_{sijk} + Q_{pqis} R_{hsjk} 
                        + Q_{pqjs} R_{hisk} + Q_{pqks} R_{hijs}\right)
$$
vanish identically. We omit the long but standard computations. 

In general, the holomorphic pseudosymmetry is essential in the sense that the structure is not semisymmetric ($R\cdot R\not=0$). For instance, the component  
$$
  (R\cdot R)_{1(2m+1)122(2m+1)} = th^2h'(h+th')
$$
is non-zero for a suitably chosen function $h$.

For the components of the Ricci curvature tensor $S$, we have
\begin{eqnarray*}
  && S_{\alpha\alpha} = S_{\alpha'\alpha'} = -2((m+2)h^2+thh'),\\
  && S_{(2m+1)(2m+1)} = S_{(2m+2)(2m+2)} = -2(m+2)h^2-(2m+7)thh'-t^2(h'^{\,2}+hh''), 
\end{eqnarray*}
$S_{ij}=0$ otherwise, and for the scalar curvature $r$, 
$$
  r =\null-4(m+1)(m+2)h^2 - 2(4m+7)thh' - 2t^2(h'^{\,2}+hh'').
$$
Hence the scalar curvature is non-constant in general. 

However, in the above way, we can obtain non-compact holomorphically pseudosymmetric K\"ahler manifolds with constant scalar curvature $r$ and $f\geqslant0$. Indeed, if we suppose 
$$
  h(t)=\frac{1}{t^{2+m}}\sqrt{a+bt^2+ct^{4+2m}},
$$
where $a,b,c$ are certain constants such that $a+bt^2+ct^{4+2m}>0$ on a certain interval $(A,B)$, $B>A>0$, then we find 
\begin{eqnarray*}
  r     &=& -4c(m+1)(m+2)=\mathop{\rm const.}, \\
  f(t)  &=& \frac{1}{t^{4+2m}}\big(a(m+1)+bmt^2-ct^{4+2m}\big).
\end{eqnarray*}
To be sure more concrete examples, 
\begin{enumerate}
\item[(i)] if $a=c=0$ and $b=1$, then $r=0$ and $f(t)=\dfrac{m}{t^{2m+2}}>0$;
\item[(ii)] if $a=1$ and $b=c=0$, then $r=0$ and $f(t)=\dfrac{m+1}{t^{2m+4}}>0$;
\item[(iii)] if $a>0$, $b=0$, $c=-a$ (here, $(A,B)=(0,1)$), then $r=4a(m+1)(m+2)$ and $f(t)=a\Big(1+\dfrac{m+1}{t^{2m+4}}\Big)>0$. 
\end{enumerate}
One can easily note that the structures are not semisymmetric, and in the cases (ii) and (iii), they are Einstein.

\bigskip
{\it Acknowledgments.} The author would like to thank the referee for his valuable remarks which improved the paper.


\bigskip
\bigskip
\noindent
Institute of Mathematics and Computer Science\\
Wroc{\l}aw University of Technology\\
Wybrze\.ze Wyspia\'nskiego 27\\
50-370 Wroc{\l}aw\\
Poland\\
E-mail: {\tt zbigniew.olszak@pwr.wroc.pl}

\end{document}